\documentclass[12pt]{article}
\usepackage{amsmath,amsthm,amsfonts,amssymb}
\title{Dual random fragmentation and coagulation and an application
to the genealogy of Yule processes}
\def\proof{\noindent{\bf Proof:}\hskip10pt}        
\def\QED{\hfill $\Box$} 
\font\tenmath=msbm10 scaled 1200
\font\sevenmath=msbm7 scaled 1200
\font\fivemath=msbm5 scaled 1200
\newfam\mathfam \textfont\mathfam=\tenmath
\scriptfont\mathfam=\sevenmath \scriptscriptfont\mathfam=\fivemath
\def\math{\fam\mathfam}
\def \\ { \cr }
\def\R{{\math R}}
\def\wt{\widetilde}
\def\N{{\math N}}
\def\E{{\math E}}
\def\P{{\math P}}

\def \e{{\rm e}}

\newtheorem{theorem}{Theorem}
\newtheorem{proposition}[theorem]{Proposition}
\newtheorem{lemma}[theorem]{Lemma}
\newtheorem{corollary}[theorem]{Corollary}

\begin{document}
\author{
Jean Bertoin \thanks{Laboratoire de Probabilit\'es et Mod\` eles Al\'eatoires
and Institut universitaire de France,
Universit\'e Pierre et Marie Curie, 175, rue du Chevaleret,
F-75013 Paris, France.}
\hspace{.2cm}
and 
\hspace{.2cm}
Christina Goldschmidt \thanks{Laboratoire de Probabilit\'es et Mod\` eles
Al\'eatoires,
Universit\'e Pierre et Marie Curie, 175, rue du Chevaleret,
F-75013 Paris, France.}}

\date{\today}
\maketitle

\abstract{The purpose of this work is to describe a duality between
a fragmentation associated to certain Dirichlet distributions and a
natural random coagulation. The dual fragmentation and coalescent
chains arising in this setting appear in the description of the
genealogy of Yule processes.}

\section{Introduction}

At a naive level, fragmentation and coagulation are inverse phenomena,
in that a simple time-reversal changes one into the other. However,
stochastic models for fragmentation and coalescence usually impose
strong hypotheses on the dynamics of the processes, such as the
branching property for fragmentation (distinct fragments evolve
independently as time passes), and these requirements do not tend to
be compatible with time-reversal. Thus, in general, the time-reversal of
a coalescent process is not a fragmentation process.

Nonetheless, there are a few special cases in which time-reversal
does transform a coalescent process into a fragmentation process.
Probably the most important example was discovered by Pitman
\cite{Pi}; it is related to the so-called cascades of Ruelle and the
Bolthausen-Sznitman coalescent \cite{BS}, and also has a natural
interpretation in terms of the genealogy of a remarkable branching
process considered by Neveu, see \cite{BLG} and \cite{BP}.

The first purpose of this note is to point out other simple instances
of such duality, which rely on certain Dirichlet and Poisson-Dirichlet
distributions.  Then, in the second part, we shall show that these examples
are related to the genealogy of Yule processes.

\section{Dual fragmentation and coagulation}

\subsection{Some notation}

For every integer $n\geq1$, we consider the simplex
$$\Delta_n\,:=\,\left\{x=(x_1,\ldots,x_{n+1}): x_i\geq0\hbox{ for every
}i=1,\ldots,n+1\hbox{ and }\sum_{i=1}^{n+1}x_i=1\right\}\,.$$
It will also be convenient to agree that $\Delta_0:=\{1\}$.
We shall often refer to the coordinates $x_1,\ldots, x_{n+1}$
of points $x$ in $\Delta_n$ as {\it masses}.

We recall that the $n$-dimensional Dirichlet distribution with
parameter $(\alpha_1,\ldots,\alpha_{n+1})$ is
the probability measure on the simplex $\Delta_{n}$ with
density
$$\frac{\Gamma(\alpha_1+\cdots+\alpha_{n+1})}
{\Gamma(\alpha_1)\cdots \Gamma(\alpha_{n+1})}x_1^{\alpha_1-1}\cdots
x_{n+1}^{\alpha_{n+1}-1}\,.$$
The special case when
$\alpha_1=\ldots=\alpha_{n+1}:=\alpha\in\ ]0,\infty[$ will have an important
role in this work; it will be convenient to write
${\rm Dir}_n(\alpha)$ for this distribution. We recall the following
well-known construction: let $\gamma_1, \ldots, \gamma_{n+1}$ be
i.i.d.\ gamma variables with  parameters
$(\alpha,c)$. Set
$\bar\gamma=\gamma_1+\cdots+\gamma_{n+1}$, so that $\bar\gamma$ has a gamma
distribution with parameters $(\alpha(n+1),c)$. 
 Then the $(n+1)$-tuple
$$\left({\gamma_1/\bar \gamma},\ldots,{\gamma_{n+1} /\bar\gamma}\right)$$
has the distribution ${\rm Dir}_n(\alpha)$ and is independent of
$\bar\gamma$.

We also define the (ranked) infinite simplex
$$\Delta_\infty\,:=\,\left\{x=(x_1,\ldots): x_1\geq x_2\geq\ldots
\geq0 \hbox{ and }\sum_{i=1}^{\infty}x_i=1\right\}\,$$
and recall that 
the Poisson-Dirichlet distribution with parameter $\theta>0$, 
which will be denoted by ${\rm PD}(\theta)$ in the sequel,
is the law of the random sequence
$$\xi:=\left(\frac{a_1}{\sum_{i=1}^{\infty}a_i},
\frac{a_2}{\sum_{i=1}^{\infty}a_i}, \ldots\right),$$ where $a_1\geq
a_2\geq \ldots > 0$ are the atoms of a Poisson random measure on
$]0,\infty[$ with intensity $\theta y^{-1} \e^{-y}dy$.  We also recall
that $\xi$ is independent of $\sum_{i=1}^{\infty}a_i$, and that the latter
has the gamma distribution with parameters $(\theta,1)$.  By the
celebrated L\'evy-It\^o decomposition of subordinators, we may also
rephrase this construction as follows: if $\gamma=(\gamma(t), t\geq0)$
is a standard gamma process and, for each fixed $\theta>0$,
$\delta_1\geq \delta_2\geq \ldots$ denotes the sequence of sizes of
the jumps of $\gamma$ on the time interval $[0,\theta]$, then
$$\left(\frac{\delta_1}{\gamma(\theta)},\frac{\delta_2}{\gamma(\theta)},
\ldots\right)$$
has the ${\rm PD}(\theta)$ distribution and is independent of
$\gamma(\theta)$.

\subsection{Two dual random transformations}
We now define two random transformations:
$${\rm Frag}_k: \Delta_n\to \Delta_{n+k} \quad \hbox{and}\quad 
{\rm Coag}_k: \Delta_{n+k}\to\Delta_n\,,$$ 
where $k,n$ are integers.  

First,
we fix $x=(x_1,\ldots,x_{n+1})\in\Delta_n$ and 
 pick an index $I\in\{1,\ldots,n+1\}$ at random according to
the distribution
$$\P(I=i)\,=\,x_i\,,\qquad i=1,\ldots,n+1\,,$$ 
so that $x_I$ is a
size-biased pick from the sequence $x$. Let
$\eta=(\eta_1,\ldots, \eta_{k+1})$ be a random variable with values in
$\Delta_k$ which is distributed according to ${\rm Dir}_k(1/k)$ and
independent of $I$.  Then we split the $I$th mass of $x$ according to
$\eta$ and we obtain a random variable in $\Delta_{n+k}$:
$${\rm Frag}_k(x)\,:=\,\left(x_1,\ldots,x_{I-1},
x_I \eta_1,\ldots,
x_I \eta_{k+1}, x_{I+1},\ldots,x_{n+1}\right)\,.$$

Second, we fix $x=(x_1,\ldots,x_{n+k+1}) \in\Delta_{n+k}$ and pick an
index $J\in\{1,\ldots,n+1\}$ uniformly at random.  We merge the $k+1$
masses $x_{J},x_{J+1}\ldots,x_{J+k}$ to form a single mass
$\sum_{i=J}^{J+k}x_{i}$ and leave the other masses unchanged.  We
obtain a random variable in $\Delta_{n}$:
$${\rm Coag}_{k}(x)\,=\,\left(x_1,\ldots,x_{J-1}, 
\sum_{i=J}^{J+k}x_{i}, x_{J+k+1},\ldots,x_{n+k+1}\right)\,.$$

\noindent{\bf Remark.}  Consider the following alternative random
coagulation of $x=(x_1,\ldots,x_{n+k+1}) \in \Delta_{n+k}$. Pick $k+1$
indices $i_1,\ldots,i_{k+1}$ from $\{1,\ldots,n+k+1\}$ uniformly at
random without replacement, merge the masses $x_{i_1},
\ldots,x_{i_{k+1}}$, leave the other masses unchanged and let
$\widetilde {\rm Coag}_{k}(x)$ be the sequence obtained by ranking the
resulting masses in decreasing order.  Write also ${\rm
Coag}^{\downarrow}_k(x)$ for the sequence ${\rm Coag}_k(x)$
re-arranged in decreasing order.  Then if $\xi$ is exchangeable the
pairs $(\xi, {\rm Coag}^{\downarrow}_{k}(\xi))$ and $(\xi, \widetilde
{\rm Coag}_{k}(\xi))$ have the same distribution.  This remark applies
in particular to the case when $\xi$ has the law ${\rm
Dir}_{n+k}(1/k)$, and can thus be combined with forthcoming
Proposition \ref{P1}.

The starting point of this work lies in the observation of a simple
relation of duality which links these two random transformations
via Dirichlet laws. 

\begin{proposition}\label{P1} Let $k, n\geq 1$ be two integers, and
$\xi, \xi'$ two random variables with values in $\Delta_n$ and
$\Delta_{n+k}$, respectively. The following assertions are then
equivalent:

(i) $\xi$ has the law ${\rm Dir}_n(1/k)$ and, conditionally on
$\xi$, $\xi'$ is distributed as ${\rm Frag}_k(\xi)$.

(ii) $\xi'$ has the law ${\rm Dir}_{n+k}(1/k)$ and, conditionally on
$\xi'$, $\xi$ is distributed as ${\rm
Coag}_k(\xi')$.
\end{proposition}

It has been observed by Kingman \cite{King} that
for $k=1$, if $\xi'$ is uniformly distributed on the simplex 
$\Delta_{n+1}$ (i.e.\ has the law ${\rm Dir}_{n+1}(1)$), then ${\rm
Coag}_1(\xi')$ is uniformly distributed on $\Delta_n$. Clearly, this
agrees with our statement.

\vskip 2mm

\proof Let $\gamma_1, \gamma_2, \ldots, \gamma_{n+1}$ be independent
Gamma($1/k,1$) random variables and set
$$\bar\gamma\,=\,\sum_{i=1}^{n+1} \gamma_i\quad \hbox{and}\quad
\xi\,=\,\left(\frac{\gamma_1}{\bar\gamma},\ldots,
\frac{\gamma_{n+1}}{\bar\gamma}\right)\,,$$ so that $\xi$ has law
${\rm Dir}_n(1/k)$ and is independent of $\bar \gamma$.  Suppose that
$\eta$ is a ${\rm Dir}_k(1/k)$ random variable which is independent of
the $\gamma_i$'s, and let $\Phi: \R^{n+k+1} \rightarrow \R$ be a
bounded measurable function.  Let $I$ be an index picked at random from
$\{1,\ldots,n+1\}$ according to the conditional distribution
$$\P(I=i\, |\, \gamma_1, \ldots,
\gamma_{n+1})\,=\,\gamma_i/\bar\gamma\,,\qquad i=1,\ldots,n+1\,,$$ and
denote by ${\rm Frag}_k(\xi)$ the random sequence obtained from $\xi$
after the fragmentation of its $I$th mass according to $\eta$. We have
\[
\E \left(\Phi({\rm Frag}_k(\xi)), I=i \right)
= \E \left[  \frac{\gamma_i}{
\bar \gamma} \Phi \left( (\gamma_l/\bar \gamma)_{l < i},\gamma_i
\eta/\bar \gamma, (\gamma_l/\bar \gamma)_{l > i} \right) \right].
\]
Now, using the independence of $\bar\gamma$ and $\xi$ and the fact that
$\bar\gamma$ has the law Gamma($(n+1)/k, 1$), we see that the last
expression is equal to
\begin{eqnarray*}
 & &\frac {k}{n+1} \E \left[ \gamma_i \Phi \left( (\gamma_l/\bar \gamma)_{l < i},\gamma_i
\eta/\bar \gamma, (\gamma_l/\bar \gamma)_{l > i} \right) \right] \\
&= &\frac{k}{n+1} \E \int_0^{\infty} x \Phi \left(
\frac{(\gamma_l)_{l<i}}{x +
\sum_{j\neq i} \gamma_j}, \frac{x \eta}{x +
\sum_{j\neq i} \gamma_j}, \frac{(\gamma_l)_{l>i}}{x +
\sum_{j\neq i} \gamma_j} \right) \frac{1}{\Gamma(1/k)} x^{1/k - 1}
\e^{-x} dx \\
& = & \frac{1}{n+1}\E \left[\Phi \left(\frac{(\gamma_l)_{l<i}}{\gamma' +
\sum_{j\neq i} \gamma_j}, \frac{\gamma' \eta}{\gamma' +
\sum_{j\neq i} \gamma_j}, \frac{(\gamma_l)_{l>i}}{\gamma' +
\sum_{j\neq i} \gamma_j} \right) \right]
\end{eqnarray*}
where $\gamma' \sim {\rm Gamma}((k+1)/k, 1)$, independently of $\eta$
and $(\gamma_j)_{j \neq i}$.  But then $\gamma'
\eta$ is a collection of $k+1$ independent Gamma($1/k,1$) random
variables, so ${\rm Frag}_k(\xi) $ has the law ${\rm Dir}_{n+k}(1/k)$
and is independent of the random index $I$ which is uniformly distributed
on $\{1,\ldots,n+1\}$.
Since we can recover $\xi$ from ${\rm Frag}_k(\xi)$ and $I$ by an obvious
coagulation, this completes the proof. 
\QED

\vskip 2mm

Next we turn our attention to the infinite ranked simplex
and define two random transformations, ${\rm Frag}_{\infty}:
\Delta_{\infty}
\to\Delta_{\infty}$ and ${\rm Coag}_a: \Delta_{\infty}\to\Delta_{\infty}$,
where $a\in [0,1]$ is some parameter. The fragmentation transformation
on the infinite simplex simply mimics that on the finite simplex;
in this direction, recall that the Poisson-Dirichlet ${\rm PD}(1)$
arises as the weak limit as $k\to\infty$ of sequence of ${\rm
Dir}_k(1/k)$ variables after obvious re-ordering.  More
precisely, given $x=(x_1,\ldots)\in\Delta_\infty$, we pick a mass 
$x_I$ at random by size-biased sampling and split $x_I$ using an
independent variable $\eta=(\eta_1,\ldots)$ with law ${\rm PD}(1)$. In
other words,
${\rm Frag}_{\infty}(x)$ is the ranked sequence with unordered terms
$x_1,\ldots,x_{I-1}, x_I\eta_1, x_I\eta_2,\ldots,x_{I+1},\ldots$. 

Next, consider a sequence $U_1,U_2,\ldots$ of i.i.d.\ uniform random
variables and $a\in[0,1]$.  Starting again from some fixed $x\in
\Delta_{\infty}$, we merge the masses $x_i$ for which $U_i\leq a$ into a
single mass and leave the others unchanged. We denote by ${\rm
Coag}_a(x)$ the random sequence obtained by putting the resulting masses
in decreasing order. We then have the following analogue of Proposition
\ref{P1}, which is reminiscent of Corollary 13 of Pitman \cite{Pi}.

\begin{proposition}\label{P2} Let $\xi, \xi'$ be two random variables 
with values in 
$\Delta_{\infty}$. For every $\theta>0$, the following assertions are
equivalent:

(i) $\xi$ has the law ${\rm PD}(\theta)$ and, conditionally on
$\xi$, $\xi'$ is distributed as ${\rm
Frag}_{\infty}(\xi)$.

(ii) $\xi'$ has the law ${\rm PD}(\theta+1)$ and, conditionally on
$\xi'$, $\xi$ is distributed as ${\rm
Coag}_{1/(\theta+1)}(\xi')$.

\end{proposition}

\proof Let $\gamma = (\gamma(t), t \geq 0)$ be a standard gamma
process and set 
$$D_t = \gamma((\theta + 1)t)/ \gamma(\theta + 1),$$ for
$0 \leq t \leq 1$, so that $(D_t, 0 \leq t \leq 1)$ is a Dirichlet
process of parameter $\theta + 1$.  (The vector of ordered jumps of
this Dirichlet process has the ${\rm PD}(\theta + 1)$ distribution.)
Consider the following alternative way of thinking of the random
coagulation operator ${\rm
Coag}_{1/(\theta+1)}$: pick a point $V$ uniformly in $[0,1]$ and define a
new process $(D'_t, 0 \leq t \leq 1)$ by
\[
D'_t = \begin{cases} 
	D_{\theta t/(\theta + 1)} & \text{ if $t < V$} \\
	D_{(1 + \theta t)/(\theta + 1)} & \text{ if $t \geq V$}.
       \end{cases}
\]
As the times of the jumps of $D$ are uniformly distributed on $[0,1]$,
this picks a proportion $1/(\theta + 1)$ of them and coalesces them into a
single jump (say $\beta^* = D_{(1 + \theta V)/(\theta + 1)}-
D_{\theta V/(\theta + 1)}$) at $V$.  Let $\beta_1 \geq \beta_2 \geq
\ldots > 0$ be the sequence of other jumps of $D'$ and $U_1, U_2,
\ldots$ the corresponding jump times.  Let $\beta_1' \geq \beta_2'
\geq \ldots > 0$ be the sequence of jumps of $D$ in the interval
$[\theta V/(\theta + 1), (1 + \theta V)/(\theta + 1)]$, so that $\beta^*
= \sum_{i=1}^{\infty} \beta'_i$.  We wish to show that $D'$ is a
Dirichlet process with parameter $\theta$, so that the vector
$(\beta^*, \beta_1, \beta_2, \ldots)$ of its jumps (re-arranged in the
decreasing order) has the ${\rm PD}(\theta)$ distribution.  We will
also show that the mass $\beta^*$
resulting from the coalescence constitutes a size-biased pick from this
vector.

Let
\begin{align*}
\gamma^1(t) & = \begin{cases}
		\gamma(t) & \text{ if $t < V \theta$} \\
		\gamma(t+1) - (\gamma(V \theta + 1) - \gamma(V
		\theta)) & \text{ if $V \theta \leq t \leq \theta$}
	       \end{cases} \\
\gamma^2(t) & = \gamma(V \theta + t) - \gamma(V \theta) \qquad \text{for $0
		\leq t \leq 1$}.
\end{align*}
Then $\gamma^1$ and $\gamma^2$ are independent processes with
$\gamma^1 \stackrel{d}{=} (\gamma(t), 0 \leq t \leq \theta)$ and
$\gamma^2 \stackrel{d}{=} (\gamma(t), 0 \leq t \leq 1)$,
independently of $V$.  Write $\delta_1 \geq \delta_2 \geq \ldots$ for
the ordered sequence of jumps of $\gamma^1$ and $T_1, T_2, \ldots$ for
the corresponding times of these jumps.  Write $\delta'_1 \geq
\delta'_2 \geq \ldots$ for the ordered sequence of jumps of
$\gamma^2$.  Then
\renewcommand{\labelenumi}{(\roman{enumi})}
\begin{enumerate}
\item $U_1 = T_1/\theta, U_2 = T_2/\theta, \ldots$ are i.i.d.\ ${\rm
U}[0,1]$,
\item $\beta^* = \gamma^2(1)/ \gamma(1 + \theta)$ and so has a ${\rm
Beta}(\theta, 1)$ distribution,
\item $\frac{1}{\beta^*} (\beta_1', \beta_2', \ldots) =
\frac{1}{\gamma^2(1)}(\delta_1', \delta_2', \ldots)$ and so has the
${\rm PD}(1)$ distribution,
\item $\frac{1}{1 - \beta^*} (\beta_1, \beta_2, \ldots) =
\frac{1}{\gamma^1(\theta)}(\delta_1, \delta_2, \ldots)$ and so has
the ${\rm PD}(\theta)$ distribution.
\end{enumerate} 
Furthermore, the random variables in (i) to (iv) above are
independent.  The fact that $\beta^*$ is a size-biased pick from
$(\beta^*, \beta_1, \beta_2, \ldots)$ and the ${\rm PD}(\theta)$
distribution of the latter follow from (i) and (iii) and the
stick-breaking scheme (see, for instance, Definition 1 in Pitman and
Yor \cite{PY}).  That $D'$ is a Dirichlet process of parameter
$\theta$ then follows from (iv) and the independence.

The coagulation operator used here can be re-phrased as follows:
starting with $x \in \Delta_{\infty}$, take a sequence $V, V_1, V_2,
\ldots$ of i.i.d.\ ${\rm U}[0,1]$ random variables, merge the masses
$x_i$ for which $V_i \in [\theta V/(\theta + 1), (1 + \theta
V)/(\theta + 1)]$, leave the other masses unchanged and, finally, rank
the resulting sequence in decreasing order.  Call this operator
$\wt{\rm Coag}_{1/(\theta + 1)}$.  Then it is clear that whenever
$\xi'$ is a random exchangeable sequence in $\Delta_{\infty}$, $(\xi',
{\rm Coag}_{1/(\theta + 1)}(\xi'))$ and $(\xi', \wt{\rm
Coag}_{1/(\theta + 1)}(\xi'))$ have the same distribution.  Our claim
follows now readily from these results.  \QED

\vskip 2mm \noindent
{\bf Remark.} It may be interesting to check 
Proposition~\ref{P2} as follows. Consider Poisson random
measure $M$ on
$(0, \infty)$ with intensity $\theta x^{-1} \e^{-x} dx$.  Let $a_1,
a_2, \ldots$ be the atoms of $M$ in decreasing order, so that 
\[
\left(\frac{a_1}{\sum_{j=1}^{\infty} a_j},
\frac{a_2}{\sum_{j=1}^{\infty} a_j}, \ldots \right)
\]
has distribution PD($\theta$), independently of $\sum_{j=1}^{\infty}
a_j$.  Let $\eta \sim {\rm PD}(1)$, independently of $M$ and suppose
that $\Phi: \Delta_{\infty} \rightarrow \R$ is any symmetric bounded
measurable function.  Then if $\xi \sim {\rm PD}(\theta)$, using
independence we have that
\[
\E \left[ \Phi({\rm Frag}_{\infty}(\xi)) \right]
= \frac{1}{\E \left[ \sum_{j=1}^{\infty} a_j \right]} \E \left[
\sum_{i=1}^{\infty} a_i \Phi \left( \frac{a_i
\eta}{\sum_{j=1}^{\infty} a_j}, \frac{(a_l)_{l \neq
i}}{\sum_{j=1}^{\infty} a_j} \right) \right].
\]
By the Palm formula, this is equal to
\begin{eqnarray*}
& &\frac{1}{\theta} \E \int_0^{\infty} x \Phi \left( \frac{x \eta}{x +
\sum_{j=1}^{\infty} a_j},
\frac{(a_l)_{l=1}^{\infty}}{x + \sum_{j=1}^{\infty} a_j} \right) \theta
x^{-1} \e^{-x} dx \\
&=& \E \left[ \Phi \left( \frac{a' \eta}{a' +
\sum_{j=1}^{\infty} a_j},
\frac{(a_l)_{l=1}^{\infty}}{a' + \sum_{j=1}^{\infty} a_j} \right) \right]
\end{eqnarray*}
where $a' \sim {\rm Exp}(1)$, independently of $M$ and $\eta$.  But
then $a' \eta$ has the distribution of the atoms of a Poisson random
measure with intensity $x^{-1} \e^{-x} dx$ arranged in decreasing order
and so we see that taking these atoms together with those of $M$, we
obtain a Poisson random measure of intensity $(\theta + 1) x^{-1}
\e^{-1} dx$.  Hence, ${\rm Frag}_{\infty}(\xi)$ has the law $ {\rm
PD}(\theta + 1)$.

\subsection{Dual fragmentation and coagulation chains} 
\label{subsec:chains}

The dual fragmentation and coagulation operators that were defined in
the preceding section incite us to introduce Markov fragmentation and
coagulation chains in duality by time-reversal. Specifically, we
consider for each integer $k\geq1$ a chain
$$X^{(k)}(0),X^{(k)}(1),X^{(k)}(2),\ldots\,,$$ 
where $X^{(k)}(n)$ is a random variable with values in $\Delta_{nk}$
(in particular $X^{(k)}(0)=1$), and the conditional distribution of
$X^{(k)}(n+1)$ given $X^{(k)}(n)=x$ is the law of ${\rm
Frag}_k(x)$. We deduce from Proposition \ref{P1} by induction that for
each $n$, $X^{(k)}(n)$ has the distribution ${\rm Dir}_{nk}(1/k)$.
The time-reversed coagulation chain
$$\ldots,X^{(k)}(n+1), X^{(k)}(n), \ldots, X^{(k)}(1), X^{(k)}(0)$$ is
also Markov; more precisely, the conditional distribution of
$X^{(k)}(n)$ given $X^{(k)}(n+1)=x$ is the law of ${\rm Coag}_k(x)$.
Note that for $k=1$, this has the distribution of the jump chain of
Kingman's coalescent \cite{King}.

Analogously, for $k=\infty$, we can define a Markov
fragmentation chain on $\Delta_{\infty}$,
$$X^{(\infty)}(0),X^{(\infty)}(1),X^{(\infty)}(2),\ldots\,,$$ such
that the conditional distribution of $X^{(\infty)}(n+1)$ given
$X^{(\infty)}(n)=x$ is the law of ${\rm Frag}_{\infty}(x)$.  We deduce
by induction from Proposition \ref{P2} that for every $\theta>0$, if
the distribution of the initial state $X^{(\infty)}(0)$ is ${\rm
PD}(\theta)$ then, for every integer $n$, $X^{(\infty)}(n)$ has the
distribution ${\rm PD}(\theta+n)$. Moreover, in this case, the
time-reversed coagulation chain
$$\ldots,X^{(\infty)}(n+1),
X^{(\infty)}(n), \ldots, X^{(\infty)}(1), X^{(\infty)}(0)$$ is also
Markov; more precisely, the conditional distribution of
$X^{(\infty)}(n)$ given $X^{(\infty)}(n+1)=x$ is the law of 
${\rm Coag}_{1/(n+1+\theta)}(x)$.

\vskip 2mm
\noindent {\bf Remarks.} (a) Recall that the parameter $\theta$ can be
recovered from a sample $\xi$ of a ${\rm PD}(\theta)$ random variable as
follows:
$$\theta\,=\,\lim_{\varepsilon\to0+}\frac{1}{\log 1/\varepsilon}
\max\left\{n: \xi_n>\varepsilon\right\}\,.
$$
This shows that the above description for the reversed coagulation chain
is indeed Markovian.

(b) There is simple representation for the $k = \infty$
fragmentation chain in terms of compound bridges with exchangeable
increments which is inspired by
\cite{BLG2}.  Let
$U_0,U_1,\ldots$ be a sequence of independent uniform variables on $[0,1]$.
For each $n$, we consider the elementary bridge $b_n: [0,1]\to[0,1]$
defined by
$$b_n(t)\,=\,\frac{n}{n+1}t + \frac{1}{n+1}{\bf 1}_{\{t>U_n\}}\,,\qquad
t\in [0,1]\,.$$
Then is is easy to check that for every $n\in\N$,
the sequence $b_n\circ b_{n+1}\circ \cdots\circ b_{n+i}$ converges
pointwise almost surely as $i\to\infty$ to a bridge with exchangeable
increments
$B_n$ which has no drift and infinitely many jumps a.s. If we write
$\beta_n\in\Delta_n$ for the sequence of the sizes of the jumps of
$B_n$ ranked the decreasing order, then the chain $(\beta_n, n\in\N)$
has the same law as $X^{(\infty)}$. We refer to \cite{BLG2} for the
necessary technical background.

\section{The genealogy of Yule processes}
We shall now show that the dual fragmentation and
coagulation chains which we introduced in the preceding section
are naturally connected to the genealogy of Yule processes.

\subsection{Discrete setting}
For every integer $k\geq1$, we write $Y^{(k)}=\left(Y^{(k)}_t,
t\geq0\right)$ for the Yule process started from $Y_0^{(k)} = 1$:
$Y^{(k)}_t$ gives the number of individuals alive a time $t$ in a
branching process in which each individual lives for an exponential
time of parameter $1$ and gives birth at its death to $k+1$ children,
which then evolve independently of one another according to the same
rules as their parent. We agree to label each child of an individual
by an integer in $\{1,\ldots,k+1\}$, which allows us to order
individuals at any generation in a consistent way: given two distinct
individuals, we may consider their most recent common
ancestor. Plainly, two different children of this ancestor are
ancestors of exactly one of these two individuals, and the labelling
of the children of the most recent common ancestor induces the order
of the individuals.

\begin{lemma}\label{L1}
The process $\left(\exp(-kt)Y^{(k)}_t, t\geq0\right)$ is a uniformly
integrable martingale and its limit $W^{(k)}$ has the {\rm
Gamma}$(1/k,1/k)$ distribution.
\end{lemma}

\proof A similar limit result is stated in Athreya \& Ney \cite{AN} on
page 130; however, the limiting distribution given there is incorrect
and so we shall provide here a detailed proof.  The martingale
property is classical, so we focus on the distribution of the limit
$W^{(k)}$. Define
$$\Phi_t(s):=\E\left(s^{Y^{(k)}_t}\right)\,.$$
The backward equation implies that
$$\frac{\partial}{\partial t}\Phi_t(s)\,=\,\Phi^{k+1}_t(s)-\Phi_t(s)
\quad ,\quad \Phi_0(s)=s\,.$$
This equation has solution
$$\Phi_t(s)\,=\,s
\e^{-t}\left(1-\left(1-\e^{-kt}\right)s^k\right)^{-1/k}\,.$$
Hence, for $\theta<0$, 
\begin{eqnarray*}
\E\left(\exp\left(\theta\e^{-kt}Y^{(k)}_t\right)\right)
\,&=&\,\exp\left(\theta \e^{-kt}\right)\e^{-t}
\left(1-\left(1-\e^{-kt}\right)
\exp\left(\theta k \e^{-kt}\right)\right)^{-1/k}\\
\,&=&\,\left[
\e^{kt}\exp\left(-\theta k \e^{-kt}\right)-\e^{kt}+1\right]^{-1/k}\,,
\end{eqnarray*}
and when $t\to\infty$, this quantity converges to
$$\left(1-k\theta\right)^{-1/k}\,=\,\left(\frac{1/k}{1/k
-\theta}\right)^{1/k},$$
which is the moment generating function of a gamma random variable with
parameters $(1/k,1/k)$. \QED

\vskip 2mm

We think of $W^{(k)}$ as the size of the terminal population.
For every $t\geq0$, by application of the branching property at time $t$,
we may decompose the terminal population into sub-populations having the
same ancestor at time $t$. Specifically,
$$W^{(k)}\,=\,\sum_{i=1}^{Y^{(k)}_t}W^{(k)}_i(t)\,,$$ where
$W^{(k)}_i(t)$ is the size of the terminal sub-population descending
from the $i$th individual in the population at time $t$. Observe that
conditionally on $Y^{(k)}_t$, the variables $W^{(k)}_i(t)$ are
independent and all have the same law as $\e^{-kt}W^{(k)}$.

Finally, we define the genealogical process $G^{(k)}=\left(G^{(k)}(t),
t\geq0\right)$ associated to $Y^{(k)}$ by
$$G^{(k)}(t)\,=\,\left(W^{(k)}_1(t), \ldots,
W^{(k)}_{Y^{(k)}_{t}}(t)\right)\,.$$ 

The genealogical structure of the Yule process can be described in
terms of the fragmentation chain $X^{(k)}$ of
Section~\ref{subsec:chains} as follows.

\begin{theorem}\label{T1}  Let 
$N=(N_t, t\geq0)$ be a standard Poisson process which is
independent of the chain $X^{(k)}$. Then for each $w>0$, the compound
chain 
$$\left(wX^{(k)}( N_{wt}), t\geq0\right)$$
has the same law as the time-changed process 
$$\left(G^{(k)}\left(\frac{1}{k}
\log(1 + kt)\right), t\geq0\right)$$ conditioned on
$W^{(k)}=w$.
\end{theorem}
\vskip 2mm
\noindent
{\bf Remark.} Theorem I of Kendall~\cite{Ken} states that given $W^{(1)}$,
$\left(Y^{(1)}_{\log(1 + t/W^{(1)})}, t \geq 0\right)$ is a Poisson
process with unit parameter.  This is clearly an aspect of
Theorem~\ref{T1}.  Moreover, on page 130 of Athreya \& Ney~\cite{AN}, it
is suggested that no generalization of Kendall's result to a more
general continuous-time Markov branching process is known;
Theorem~\ref{T1} constitutes a small such generalization.

\vskip 2mm

\proof Set $\tau(t):=\frac{1}{k}
\log(1 + kt)$ and 
let $T$ be the time of the first birth in the Yule process, which is also
the time of the first dislocation of $G^{(k)}$.  The $k + 1$
fragments of
$G^{(k)}(T)$ can be written as $\e^{-kT}Z_1, \ldots, \e^{-kT}Z_{k+1}$
where, by the branching property, $Z_1, \ldots, Z_{k+1}$ are i.i.d.\ ${\rm
Gamma}(1/k,1/k)$ random variables, independent of $T$ which is
Exp($1$).  Define a change of variables by
\begin{align*}
S & = \tau^{-1}(T) = (\e^{kT}-1)/k \\
U_1 & = \e^{-kT} Z_1, \quad \ldots, \quad U_k =
\e^{-kT} Z_k, \quad W = \e^{-kT} (Z_1 + \cdots + Z_{k+1}).
\end{align*}
It is straightforward to see that the joint density of $(T, Z_1,
\ldots, Z_{k+1})$ is
\begin{align*}
& f(t, z_1, \ldots, z_{k+1}) \\
&  = \e^{-t} \Gamma(1/k)^{-(k+1)}
(1/k)^{(k+1)/k} (z_1 z_2 \ldots z_{k+1})^{-(k-1)/k} \exp(-(z_1 + \cdots + z_{k+1})/k)
\end{align*}
and so the joint density of $(S, U_1, \ldots, U_k, W)$ is
\begin{align*}
g(s, u_1, \ldots, u_k, w) & = w \e^{-ws} \cdot (1/k)\Gamma(1/k)^{-k} w^{-1/k}
(u_1 u_2 \ldots u_k(w - u_1 - \cdots - u_k))^{1/k - 1} \\
& \quad \cdot (1/k)^{1/k} \Gamma(1/k)^{-1} w^{-(k-1)/k} \exp(-w/k).
\end{align*}
Hence, $W \sim {\rm Gamma}(1/k,1/k)$ (as we already knew) and,
conditional on $W = w$, we have $S \sim {\rm Exp}(w)$ and $(U_1, U_2,
\ldots, U_k, W - U_1 -  \cdots - U_k) \sim w {\rm Dir}_k(1/k)$
independently of $S$.  Thus, the first dislocation has the correct
dynamics.  But by the branching property, subsequent dislocations are
independent for different sub-populations and the total rate of
fragmentation is always $w$.  Hence result. \QED

\vskip 2mm

In the terminology of \cite{Be}, Theorem \ref{T1} states that the
time-changed genealogical process $G^{(k)}\circ \tau$ is a
self-similar fragmentation with index $1$, dislocation law ${\rm
Dir}_k(1/k)$ and erosion coefficient 0. It may be interesting to
observe that in the special case $k=1$, this result can also be
derived as follows.

Consider a real Brownian motion $B$ started from $1$ and killed when
it reaches $0$ (at time $T_0 = \inf\{t \geq 0: B_t = 0\}$).  For every
$u\in[0,1[$, let $\wt Y_u$ denote the number of excursions of $B$ away
from $1$ which go below level $u$.  Then $(Y^{(1)}_{-\log(1 - u)})_{0
\leq u < 1}$ is a version of $({\wt Y_u})_{0 \leq u < 1}$.  To see
this, let us consider the evolution of $\wt Y$.  Firstly, ${\wt Y_0} =
1$, corresponding to the single excursion below $1$ which reaches $0$.
Let $D = \sup \{t < T_0: B_t = 1\}$, the starting time of the final
excursion which hits 0, let $U = \inf_{0 \leq t \leq D} B_t$ be the
level reached by the deepest excursion below 1 before $D$ and let
$T_U$ be the time at which it is reached.  Then, by Williams' path
decomposition theorem (Theorem VII.4.9 of Revuz and Yor~\cite{RY}),
$U$ is distributed uniformly on $[0,1[$ and, conditional on $U$,
$(B_t)_{0 \leq t < T_U}$ is a Brownian motion started at $1$ and
stopped when it first hits level $U$.  By symmetry, $(B_{D - t})_{0
\leq t < D-T_U}$ is another independent Brownian motion started at $1$
and stopped when it first hits level $U$.  Thus, ${\wt Y_u}$ is equal
to 1 on $[0, U[$, ${\wt Y_U} = 2$ and $({\wt Y_{U + v}})_{0 \leq v <
1-U}$ evolves as the sum of two independent processes which are the
same as $\wt Y$ except that the times until the first jumps are now
uniform on $[0,U[$ rather than on $[0,1[$.  (This is Theorem 8 of Le
Gall~\cite{LG}, repeated here for completeness.)  Time-changing
$Y^{(1)}$ with $u \to -\log(1 - u)$ means that its exponential
inter-jump times become uniform and so we do, indeed, have
$(Y^{(1)}_{-\log(1 - u)})_{0 \leq u < 1} \stackrel{d}{=} ({\wt
Y_u})_{0 \leq u < 1}$.

A more elegant way of expressing the preceding argument is to say that
the Brownian path encodes a continuous-state branching process with
quadratic branching mechanism.  The local time at level 1, $L^1_{T_0}$,
satisfies
\begin{equation*}
L^1_{T_0} = \lim_{u \to 1-} 2 (1 - u) {\wt Y_u}.
\end{equation*}
In this context, $\frac{1}{2} L_{T_0}^1$ corresponds to the size of
the population at time 1 in the continuous-state branching process
generated by a single ancestor conditioned to have descendents up to
time $1$.  The so-called reduced tree associated with the population
at time $1$ is described up to the deterministic time-change
$u\to-\log(1-u)$ by the Yule process $Y^{(1)}$. See, for instance,
Section 2.7 in Duquesne and Le Gall \cite{DLG}, and Fleischmann and
Siegmund-Schultze \cite{FSS}.  Note that the well-known fact that
$\frac{1}{2} L_{T_0}^1$ has an exponential distribution with mean 1
(Proposition VI.4.6 of Revuz and Yor~\cite{RY}) gives another
derivation of the limiting distribution in Lemma \ref{L1}, since
$$ W^{(1)} = \lim_{t \to \infty} e^{-t} Y^{(1)}_t = \lim_{u \to 1-} (1 - u)
Y^{(1)}_{-\log(1 - u)} \stackrel{d}{=} \frac{1}{2} L_{T_0}^1.$$

It is known from excursion theory that in the scale of the local time
at level $1$, the rate of excursions of $B$ away from $1$ which reach
level $u\in\ ]0,1[$ but do not exceed $u-du$ is $(1-u)^{-2}du$. Note
that the map $s\to 1-\frac{1}{1+s}$ from $\R_+$ to $[0,1[$ has inverse
$u \to \frac{1}{1-u}-1$ and, thus, transforms Lebesgue measure on
$\R_+$ into the measure $(1-u)^{-2}du$ on $[0,1[$.  Suppose that we
split the local time at level $1$ according to the occurrence of
excursions exceeding level $u$, so that we obtain the sequence
$$\wt W(u)\,=\,\left(\wt W_1(u),\ldots,\wt W_{\wt
Y_u}(u)\right),$$ where $\wt W(u)$ is the sequence of the
increments of the local time at level $1$ on the maximal time
intervals such that at the beginning and end of each interval $B$ is at
$1$ and during the interval remains above level $u$.  Then it follows
easily that the time-changed process $\left(\wt W
\left(1-\frac{1}{1+s}\right), s \geq 0 \right)$ is a fragmentation in
which each mass, say $x$, splits at rate $x$ into $xU$ and $x(1-U)$
where $U$ is uniform.  In other words, conditionally on $\frac{1}{2}
L_{T_0}^1 =w$, the
process $\left(\wt W\left(1-\frac{1}{1+s}\right), s \geq 0 \right)$ is
distributed as the compound fragmentation chain $\left(wX^{(1)}(
N_{ws}), s\geq0\right)$, where $N$ is an independent standard Poisson
process.

Finally, the composition of the two time-changes which appear in this
analysis yields
$$s \to -\log \left( 1 - \left(1 - \frac{1}{1+s} \right) \right) =
\log(1 + s)\,, \qquad s\in \R_+\,,$$
and so we recover Theorem \ref{T1} in the special case $k=1$.
Unfortunately, it does not seem that there are similar interpretations for
$k\geq2$. 

\begin{corollary} \label{C1}
We have that 
$$\left(\frac{1}{W^{(k)}} G^{(k)}\left(\frac{1}{k} \log \left(1 +
ke^{-t}/W^{(k)}\right) \right), t \in \R\right)$$ is a time-homogeneous
Markov coalescent process which is independent of $W^{(k)}$.  For
any $n \geq 1$, given that it is in state $x \in \Delta_{nk}$, it
waits an exponential time of parameter $n$ and then jumps to a
variable distributed as ${\rm Coag}_{k}(x)$, independently of the
exponential time.
\end{corollary}

Note that the case $k=1$ of this result gives a variation of Kingman's
coalescent.  The jump-chains are identical, as we have already noted,
but here the rate of coalescence of two blocks depends on the total
number of blocks present, whereas in Kingman's coalescent it does not.

\vskip 2mm

\proof  Firstly, we note that by Theorem~\ref{T1},
$$\left(\frac{1}{W^{(k)}} G^{(k)}\left(\frac{1}{k} \log (1 + ke^{-t}/W^{(k)})
\right), t \in \R\right)$$
has the same law as
$$ \left(X^{(k)}(N_{e^{-t}}), t \in \R \right)$$ and so we will work
with the latter process instead.  The $k=1$ case is essentially
treated in \cite{Be2} and the proof proceeds in the same way here.
The jump chain clearly behaves in the correct manner and so it remains
to check that the inter-jump times are as claimed.  Let $0 \leq T_1
\leq T_2 \leq \ldots$ be the jump times of $(N_t)_{t \geq 0}$.  Then
the first instant that $X^{(k)}(N_{e^{-t}})$ has exactly $nk+1$ terms
is
$$ \inf \left\{ t \in \R : N_{e^{-t}} = n \right\} = -\log
T_{n+1}. $$
The sequence of inter-jump times is
$$ \ldots, \log T_{n+1} - \log T_n, \log T_n - \log T_{n-1}, \ldots,
\log T_2 - \log T_1$$
and it is easily shown that this is a sequence of independent
exponential random variables with parameters
$$ \ldots, n, n-1, \ldots, 1$$
respectively. \QED

\subsection{Continuous setting}

 {\it Continuous-state branching processes} (or CSBP's) were
introduced by Lamperti~\cite{lamperti, lamperti2} as limits of rescaled
branching processes. Typically, a CSBP is a
time-homogeneous Markov process with values in
$\R_+$,
$$Z\,=\,(Z(t,a),\, t\geq0 \hbox{ and } a\geq0)\,,$$ 
(where the parameter $t$ refers to time and the parameter $a$ to the
starting point i.e.\ $Z(0,a)=a$ a.s.) which fulfils the branching
property: the path-valued process $(Z(\cdot,x), x\geq0)$ has
independent and stationary increments. In particular, if $\wt
Z(\cdot,y)$ is an independent copy of $Z(\cdot,y)$, then
$Z(\cdot,x)+\wt Z(\cdot,y)$ has the law of $Z(\cdot,x+y)$. There is a
simple relation connecting CSBP's and Bochner's subordination for
subordinators which enables us to define their genealogy; we refer the
interested reader to \cite{BLG} for heuristics, detailed arguments
etc.

We call a {\it continuous state Yule process} a CSBP
$$Y\,=\,(Y(t,a),\, t\geq0 \hbox{ and } a\geq0)\,,$$ which evolves as
follows: for each $a>0$, the process $Y(\cdot,a)$ waits an exponential
time with parameter $a$ and then jumps to $a+1$.  It then
evolves independently as if it had been started in state $a+1$.  In
terms of the genealogy, the sub-population of size 1 which is born at
a jump time has a parent which is chosen uniformly at random from the
population present before the jump.  Note that this genealogy is easy
to describe in a consistent manner for different values $a$ of the
starting population.

It is immediate that for an integer starting point $a\in \N$, the
process $(Y(t,a), t\geq0)$ is a Yule process $Y^{(1)}$ with
$2$ offspring, as considered in the preceding section.  However, we stress
that its genealogy is not the same as that of $Y^{(1)}$,
as we are dealing with a continuous population in the first case and a
discrete population in the second.

We have the following analogue of Lemma \ref{L1}:

\begin{lemma}\label{L2} For every $a\geq0$, the process
$\left(\e^{-t}Y(t,a), t\geq0\right)$ is a uniformly integrable
martingale. Its limit, say $\gamma(a)$, viewed as a process in the
variable
$a$, has the same finite dimensional laws as a standard gamma process.
\end{lemma}

\proof 
For $a=1$, we see from Lemma \ref{L1} and the identity in distribution
$Y(\cdot,1){\mbox{$ \ \stackrel{\cal L}{=}$ }}Y^{(1)}(\cdot)$
that $\left(\e^{-t}Y(t,1), t\geq0\right)$ is a uniformly
integrable martingale and that its limit has the standard exponential
distribution.  The proof is easily completed by an appeal to the
branching property. \QED

\vskip 2mm

\noindent
{\bf Remark.} The limiting distribution in Lemma~\ref{L2} is
essentially a corollary of Theorem 3 of Grey~\cite{Grey}.

\vskip 2mm

Just as in the preceding section, we think of $\gamma(a)$ as the size of
the terminal population when the initial population has size $a$. 
We can express $\gamma(a)$ as 
$$\gamma(a)\,=\,\sum_{b\leq a}\delta_b\,,$$
where $\delta:=(\delta_b, b\geq0)$ is the jump process of $\gamma$,
which corresponds to decomposing
the terminal population into sub-populations
having the same ancestor at the initial time.
We write $G(0,a)$ for the sequence of the jumps of $\gamma$ on $[0,a]$,
ranked in decreasing order, and we deduce from Lemma \ref{L2}
that conditionally on $\gamma(a)=g$, $G(0,a)/g$
has distribution ${\rm PD}(a)$.

More generally, by the branching property, we can decompose the terminal
population into sub-populations having the same ancestor at any given time
$t$. This gives
$$\gamma(a)\,=\,\sum_{b\leq Y(t,a)}\e^{-t}\delta^{(t)}_b\,,$$
where $\delta^{(t)}:=(\delta^{(t)}_b, b\geq0)$ is the jump process of
a standard gamma process $\gamma^{(t)}$ which is independent of the Yule
process up to time $t$, $(Y(s,c),\, s \in [0,t] \hbox{ and }c\geq0)$.
This enables us to define for each $a>0$ the genealogical process
 associated to a Yule process
$Y(\cdot,a)$,
$$G(\cdot,a)=\left(G(t,a), t\geq0\right)\,,$$
where $e^t G(t,a)$
is the ranked sequence of the sizes of the jumps of the subordinator
$\gamma^{(t)}$ on the interval $[0,Y(t,a)]$.

An easy variation of the arguments for the proof of Theorem \ref{T1}
shows that the genealogical structure of the Yule process can be described
in terms of the fragmentation chain
$X^{(\infty)}$ of Section~\ref{subsec:chains} as follows.

\begin{theorem}\label{T2} Fix $a,g>0$ and let 
the chain $X^{(\infty)}$ have initial distribution
${\rm PD}(a)$. Introduce a standard Poisson process, $N=(N_t, t\geq0)$, 
which is independent of the chain
$X^{(\infty)}$. 
 Then the compound
chain 
$$\left(g X^{(\infty)}( N_{g t}), t\geq0\right)$$
has the same law as the time-changed process 
$$\left(G\left(
\log(1 + t), a \right), t\geq0\right)$$ conditioned on
$\gamma(a)=g$.
\end{theorem}

Likewise, the analogue of Corollary~\ref{C1} is as follows.

\begin{corollary}
Fix $a > 0$.  Then
$$\left(\frac{1}{\gamma(a)} G\left(\log(1 + e^{-t}/\gamma(a)),
a\right), t \in \R \right)$$ is a time-homogeneous Markov
coalescent process which is independent of $\gamma(a)$.  Suppose that it is in
state $x \in \Delta_{\infty}$ and recall Remark (a) of
Section~\ref{subsec:chains}.  Then if
$$\lim_{\epsilon \to 0+} \frac{1}{\log 1/\epsilon } \max\{i :
x_i > \epsilon\} = n + a,$$
the process waits an exponential time of parameter $n$ and then jumps
to a variable distributed as ${\rm Coag}_{1/(n+a)}(x)$, independently
of the exponential time.
\end{corollary}

\end{document}